\documentclass[reqno,12pt]{amsart}
\usepackage[cp1252]{inputenc}
\usepackage{amssymb}
\usepackage{amsmath}
\usepackage{amsfonts}

\newtheorem{theorem}{Th\'eor\`eme}[section]
\theoremstyle{plain}

\newtheorem{corollary}[theorem]{Corollaire}

\newtheorem{definition}[theorem]{D\'efinition}

\newtheorem{remark}[theorem]{Remarque}

\numberwithin{equation}{section}

\input{tcilatex}

\begin{document}
\title{La Propri\'{e}t\'{e} De Radon-Nikodym dans l'espace $\mathcal{L}(X,Y)$%
}
\author{Daher Mohammad}
\address{16, Square Albert Scheitzer- 77350 Le M\'{e}e sur Seine-France}
\email{m.daher@orange.fr}

\begin{abstract}
Soient $X$ un espace de Banach s\'{e}parable et $Y$ un espace de Banach a la
propri\'{e}t\'{e} de Radon-Nikodym. Dans ce travail on montre que $\mathcal{L%
}(X,Y)$ a la propri\'{e}t\'{e} de Radon-Nikodym si $\mathcal{L}(X,Y)$ est un
espace faiblement $LUR,$ ou si $\mathcal{L}(X,Y)$ est un espace $WCG.$

\textsc{Abstract}: Let $X$ be a separable Banach space and $Y$ a space which
has the Radon-Nikodym property. In this work, we show that $\mathcal{L}(X,Y)$
has the Radon-Nikodym property, if $\mathcal{L}(X,Y)$ is weakly locally
uniformly convex or if $\mathcal{L}(X,Y)$ is a weakly compactly generated
space.
\end{abstract}

\maketitle


Mots cl\'{e}s~: Interpolation des espaces de Hardy

\begin{center}
\textsc{Introduction}
\end{center}

Soit $X$ un espace de Banach. Dans \cite{Die-F} on montre que $X^{\ast }$ a
la propri\'{e}t\'{e} de Radon-Nikodym si $X^{\ast }$ est un espace
faiblement $LUR$. Soient $X$ un espace de Banach s\'{e}parable et $Y$ un
espace de Banach ayant la prori\'{e}t\'{e} de Radon-Nikodym. Dans ce
travail, nous montrons que $\mathcal{L}(X,Y)$ a la propri\'{e}t\'{e} de
Radon si $\mathcal{L}(X,Y)$ est un espace faiblement $LUR.$ Dans la suite,
nous montrons que si $\mathcal{L}(X,Y)$ est un espace $WCG$, alors $\mathcal{%
L}(X,Y)$ a la propri\'{e}t\'{e} de Radon-Nikodym, ce qui g\'{e}n\'{e}ralise
le r\'{e}sultat obtenu dans \cite{Steg}.

D\'{e}signons par $m$ la mesure de Lebesgue normalis\'{e}e sur le tore $%
\mathbb{T}$ et par $\mathbb{D}$ le disque unit\'{e} ouvert de $\mathbb{C}$.

Soit $X$ un espace de Banach. Pour tout $p\geq 1$ d\'{e}signons par $\mathbf{%
h}^{p}(\mathbb{D},X)$ l'espace des fonctions $h:\mathbb{D}\rightarrow X$
harmoniques telles que $\sup_{z\in \mathbb{D}}\left\Vert f(z)\right\Vert $%
\TEXTsymbol{<}$\infty ,$ si $p=\infty $ ou v\'{e}rifiant, si $p\in \left[
1,\infty \right[ ,$ sup$_{0\leq r<1}\dint\limits_{\mathbb{T}}\left\Vert
f(re^{it})\right\Vert ^{p}dm(t)<\infty $.

Pour tout Banach $X$, nontons $B_{X}$ la boule unit\'{e} ferm\'{e}e de $X.$

Notons $P_{z}$ le noyau de Poisson au point $z\in \mathbb{D}.$

\begin{definition}
\label{ik}Soient $X,Y$ deux espaces de Banach et $T:X\rightarrow Y$ un op%
\'{e}rateur born\'{e}. $T$ est un op\'{e}rateur de Radon-Nikodym si pour
toute fonction $f$ dans $\mathbf{h}^{\infty }(\mathbb{D},X),$ la suite $%
\left[ T(f(re^{it}))\right] _{r\in \left[ 0,1\right[ }$ converge dans $Y,$
pour presque tout $t\in \mathbb{T}.$
\end{definition}

Un espace de Banach $X$ a la propri\'{e}t\'{e} de Radon-Nikodym, si l'identit%
\'{e} de $X$ est un op\'{e}rateur de Radon-Nikodym.

D'apr\`{e}s \cite{B-D}, $X$ a la propri\'{e}t\'{e} de Radon-Nikodym, si et
seulement si, $\mathbf{h}^{p}(\mathbb{D},X)=L^{p}(\mathbb{T},X)$ pour un $%
p>1 $.

\bigskip

\begin{definition}
\label{bp}Soit $Y$ un espace de Banach. On dit que $Y$ poss\`{e}de d'une
expansion de \ l'identit\'{e} de dimension finie inconditionnelle, s'il
existe une suite d'op\'{e}rateurs born\'{e}s $(A_{n})_{n\geq 0}$ dans $%
\mathcal{L}(Y)$ du rang finis telle que pour tout $y\in Y$ la s\'{e}rie $%
\underset{n\geq 0}{\dsum }A_{n}y$ converge inconditionnellement vers $y$.
\end{definition}

\begin{theorem}
\label{edg}Supposons que $X^{\ast }$ a la propri\'{e}t\'{e} de
Radon-Nikodym, $Y$ poss\`{e}de d'une expansion de l'identit\'{e} de
dimension finie inconditionnelle et que $K(X,Y)$ ne contient pas $c_{0}$
isomorphiquement, alors $\mathcal{L}(X,Y)$ a la propri\'{e}t\'{e} de
Radon-Nikodym.
\end{theorem}

D\'{e}monstration.

D'apr\`{e}s \cite[th.6]{Kalt},  $\mathcal{L}(X,Y)=K(X,Y).$ Il suffit de
montrer que toute fonction dans $\mathbf{h}^{2}\left[ \mathbb{D},K(X,Y)%
\right] $ admet des limites radiales presque-partout.

Soit $f\in \mathbf{h}^{2}(\mathbb{D},K(X,Y)).$ Comme $Y$ poss\`{e}de d'une
expansion de l'identit\'{e} de dimension finie inconditionnelle, il existe
une suite d'op\'{e}rateurs born\'{e}s $(A_{n})_{n\geq 0}$ du rang finis
telle que $f(z)=\underset{k\geq 0}{\dsum }A_{k}(f(z)),$ $z\in \mathbb{D}.$
L'espace $X^{\ast }$ a la propri\'{e}t\'{e} de Radon-Nikodym, donc il existe
un sous-ensemble mesurable $\Omega $ de $\mathbb{T}$ de mesure \'{e}gale 
\`{a} 1 tel que pour tout $k\geq 0$, $A_{k}(fr(e^{i\theta }))\underset{%
r\rightarrow ^{-}1}{\rightarrow }\psi _{k}(\theta ),$ dans $K(X,Y),$ pour
tout $\theta \in \Omega .$

\emph{Etape 1:} Montrons que la fonction 
\begin{equation}
\text{ }U:\theta \in \Omega \rightarrow \sup \left\{ \left\Vert \underset{J}{%
\dsum }\psi _{k}(\theta )\right\Vert _{K(X,Y)};\text{ }J\text{ est un
sous-ensemble fini de }\mathbb{N}\right\} \text{ }  \label{a}
\end{equation}

est dans $L^{2}(\mathbb{T}).$

Pour presque tout $\theta \in \mathbb{T}$ on a

\begin{eqnarray*}
&&\sup \left\{ \left\Vert \underset{J}{\dsum }\psi _{k}(\theta )\right\Vert
_{K(X,Y)};\text{ }J\text{ est un sous-ensemble fini de }\mathbb{N}\right\} \\
&\leq &\sup \left\{ \left\Vert \underset{J}{\dsum }A_{k}f(re^{i\theta
})\right\Vert _{K(X,Y)};\text{ }J\text{ est un sous-ensemble fini de }%
\mathbb{N},\text{ }r\in \left[ 0,1\right[ \right\} \\
&\leq &\sup \left\{ \left\Vert \underset{J}{\dsum }A_{k}\right\Vert
_{K(Y,Y)};\text{ }J\text{ est un sous-ensemble fini de }\mathbb{N}\right\}
\times sup_{0<r<1}\left\Vert f(re^{i\theta })\right\Vert _{K(X,Y)},
\end{eqnarray*}

comme $Y$ poss\'{e}de d'une expansion de l'identit\'{e} de dimension finie
inconditionnelle on a%
\begin{equation*}
\sup \left\{ \left\Vert \underset{J}{\dsum }A_{k}\right\Vert _{K(Y,Y)};\text{
}J\text{ est un sous-ensemble fini de }\mathbb{N}\right\} <+\infty \text{ 
\cite[chap.VI,p.232]{Ch-Al}.}
\end{equation*}

D'autre part, d'apr\`{e}s l'in\'{e}galit\'{e} de Hardy-Littlewood \cite[th.26%
]{H-L} la quantit\'{e} $f^{\ast }(\theta )=sup_{0\leq r\leq 1}\left\Vert
f(re^{i\theta })\right\Vert _{K(X,Y)}$ existe pour presque tout $\theta \in 
\mathbb{T}$ et $f^{\ast }\in L^{2}(\mathbb{T}).$ Par cons\'{e}quent on a (%
\ref{a}).

\emph{Etape 2: }Montrons que pour presque tout $\theta \in \mathbb{T}$ la s%
\'{e}rie$\underset{k\geq 1}{\text{ }\dsum }\psi _{k}(\theta )$ converge dans 
$K(X,Y)$ pour presque tout $\theta \in \mathbb{T}$.

Par hypoth\`{e}se $K(X,Y)$ ne contient pas $c_{0}$ isomorphiquement et pour
presque tout $\theta \in \mathbb{T}$

\begin{equation*}
\sup \left\{ \left\Vert \underset{J}{\dsum }\psi _{k}(\theta )\right\Vert
_{K(X,Y)};\text{ }J\text{ est un sous-ensemble fini de }\mathbb{N}\right\}
<+\infty \text{ (d'apr\`{e}s (\ref{a})).}
\end{equation*}

\emph{\ }D'apr\`{e}s \cite[Prop.3]{Kalt}, la s\'{e}rie$\underset{k\geq 1}{%
\text{ }\dsum }\psi _{k}(\theta )$ converge dans $K(X,Y)$ pour presque tout $%
\theta \in \mathbb{T}$.

\emph{Etape 3: }Montrons que $f$ admet des limites radiales presque-partout.

En utilisant (\ref{a}) et le th\'{e}or\`{e}me de convergence domin\'{e}e, on
voit que%
\begin{eqnarray*}
\dint\limits_{\mathbb{T}}\left[ \underset{k\geq 0}{\dsum }\psi _{k}(t)\right]
P_{z}(t)dm(t) &=&\underset{k\geq 0}{\dsum }\left[ \dint\limits_{\mathbb{T}}%
\left[ \psi _{k}(t)\right] P_{z}(t)dm(t)\right] \\
\underset{k\geq 0}{\dsum }A_{k}(f(z)) &=&f(z),\text{ }\forall z\in \mathbb{D}%
.
\end{eqnarray*}

Il en r\'{e}sulte d'apr\`{e}s \cite[th.7.6]{Ru} que $f(re^{i\theta })%
\underset{r\rightarrow ^{-}1}{\rightarrow }\underset{k\geq 0}{\dsum }\psi
_{k}(\theta )$ pour presque tout $\theta \in \mathbb{T}.\blacksquare $

\begin{remark}
\label{vcd} Le th\'{e}or\`{e}me \ref{edg} nous montre que tout espace de
Banach qui ne contient pas $c_{0}$ isomorphiquement et qui n'a pas la propri%
\'{e}t\'{e} de Radon-Nikodym, ne poss\`{e}de pas d'une expansion de l'identit%
\'{e} de dimension finie inconditionnelle. En particulier un tel espace ne
poss\`{e}de pas d'une base inconditionnelle.
\end{remark}

\begin{definition}
\label{vc}Soient $X,Y$ deux espaces de Banach et $V:X\rightarrow Y$ un op%
\'{e}rateur born\'{e}. On dit que $V$ est un opr\'{e}rateur faiblement $LUR$%
, si pour toute suite $(x_{n})_{n\geq 1}$ born\'{e}e dans $X$ telle que $%
\frac{\left\Vert x_{n}\right\Vert ^{2}}{2}+\frac{\left\Vert x\right\Vert ^{2}%
}{2}-\left\Vert \frac{(x_{n}+x)}{2}\right\Vert ^{2}\underset{n\rightarrow
+\infty }{\rightarrow }0,$ alors $V(x_{n})\underset{n\rightarrow +\infty }{%
\rightarrow }V(x)$ faiblement dans $Y.$
\end{definition}

Dans la suite, nous fixons deux espaces de Banach $X,Y.$

\begin{theorem}
\label{oq}Soient $Z$ un espace de Banach et $V:\mathcal{L}(X,Y)\rightarrow Z$
un op\'{e}rateur faiblement $LUR.$ Supposons que $X$ soit s\'{e}parable et $Y
$ a la propri\'{e}t\'{e} de Radon-Nikodym, alors $V$ est un op\'{e}rateur de
Radon-Nikodym.
\end{theorem}

\emph{Etape 1:}\textbf{\ }Soit $f\in \mathbf{h}^{2}\left[ \mathbb{D},%
\mathcal{L}(X,Y)\right] .$ Montrons qu'il existe un sous-ensemble $\Omega $
de $\mathbb{T}$ de mesure \'{e}gale \`{a} 1 et une fonction $\phi
_{f}:\Omega \rightarrow \mathcal{L}(X,Y)$ tels que pour tout $\theta \in
\Omega $ et tout $x\in X$ $\ f(re^{i\theta })x\underset{r\rightarrow ^{-}1}{%
\rightarrow }\phi _{f}(\theta )x.$

Soit $(x_{n})_{n\geq 0}$ une suite dense dans $X.$ Consid\'{e}rons $X_{1}$
le $\mathbb{Q}\times \mathbb{Q}$ sous-espace vectoriel engendr\'{e} par $%
\left\{ x_{n};\text{ }n\in \mathbb{N}\right\} $ dans $X$. Comme $Y$ a la
propri\'{e}t\'{e} de Radon-Nikodym, il existe un sous-ensemble $\Omega $
dans $\mathbb{T}$ de mesure est \'{e}gale \`{a} 1 tel que pour tout $\theta
\in \Omega $ et tout $x\in X_{1}$ $f(re^{i\theta })(x\underset{r\rightarrow
^{-}1}{)\rightarrow }\psi _{f}(\theta )(x)$ dans $Y.$ D'autre part, d'apr%
\`{e}s l'in\'{e}galit\'{e} de Hardy-Littlewood \cite[th.26]{H-L} pour
presque tout $\theta \in \mathbb{T}$ et tout $x\in X_{1}$ on a 
\begin{equation*}
\left\Vert \psi _{f}(\theta )(x)\right\Vert _{Y}\leq sup_{0\leq
r<1}\left\Vert f(re^{i\theta })(x)\right\Vert _{Y}\leq sup_{0\leq
r<1}\left\Vert f(re^{i\theta })\right\Vert _{\mathcal{L}(X,Y)}\times
\left\Vert x\right\Vert =C(\theta )\times \left\Vert x\right\Vert .
\end{equation*}

o\`{u} $C(\theta )=sup_{0\leq r<1}\left\Vert f(re^{i\theta })\right\Vert _{%
\mathcal{L}(X,Y)}$ . Il en r\'{e}sulte que $\psi _{f}(\theta )$ se prolonge
d'une mani\`{e}re unique \`{a} un op\'{e}rateur $\phi _{f}(\theta )\in 
\mathcal{L}(X,Y)$.

Soient $x\in X$ et $x_{1}\in X_{1}.$ On a alors%
\begin{eqnarray}
\left\Vert f(re^{i\theta })(x)-\phi _{f}(\theta )(x)\right\Vert _{Y} &\leq
&\left\Vert f(re^{i\theta })(x_{1})-\psi _{f}(\theta )(x_{1})\right\Vert
_{Y}+\left\Vert f(re^{i\theta })(x)-f(re^{i\theta })(x_{1})\right\Vert _{Y}+
\notag \\
\left\Vert \psi _{f}(\theta )(x_{1})-\phi _{f}(\theta )(x)\right\Vert _{Y}
&\leq &\left\Vert f(re^{i\theta })(x_{1})-\psi _{f}(\theta
)(x_{1})\right\Vert _{Y}+  \notag \\
&&C(\theta )\times \left\Vert x-x_{1}\right\Vert _{X}  \label{k} \\
&&+\left\Vert \psi _{f}(\theta )(x_{1})-\phi _{f}(\theta )(x)\right\Vert
_{Y}.  \notag
\end{eqnarray}

L'in\'{e}galit\'{e} \ref{k} nous montre que $f(re^{i\theta })(x)\underset{%
r\rightarrow ^{-}1}{\rightarrow }\phi _{f}(\theta )(x)$ pour presque tout $%
\theta \in \mathbb{T}.$

\emph{Etape 2:}\textbf{\ }Montrons que $\left\Vert f\right\Vert _{\mathbf{h}%
^{2}(\mathbb{D},\mathcal{L}(X,Y))}^{2}=\dint\limits_{\mathbb{T}}\left\Vert
\phi _{f}(\theta )\right\Vert _{\mathcal{L}(X,Y)}^{2}dm(\theta ).$

Remarquons d'abord que la fonction $\theta \rightarrow \left\Vert \phi
_{f}(\theta )\right\Vert _{\mathcal{L}(X,Y)}$ est mesurable, car $B_{X}$ est
s\'{e}parable.

Soit $x\in B_{X}.$ D'apr\`{e}s \cite{B-D}, on a $f(z)(x)=\dint\limits_{%
\mathbb{T}}\phi _{f}(t)(x)P_{z}(t)dm(t),$ $z\in \mathbb{D},$ cela implique
que $\left\Vert f(z)(x)\right\Vert _{Y}^{{}}\leq \dint\limits_{\mathbb{T}%
}\left\Vert \phi _{f}(t)(x)\right\Vert _{Y}P_{z}(t)dm(t),$ donc $\left\Vert
f(z)\right\Vert _{\mathcal{L}(X,Y)}^{2}\leq \dint\limits_{\mathbb{T}%
}\left\Vert \phi _{f}(t)\right\Vert _{\mathcal{L}(X,Y)}^{2}P_{z}(t)dm(t).$
Par cons\'{e}quent%
\begin{eqnarray*}
sup_{0\leq r<1}\dint\limits_{\mathbb{T}}\left\Vert f(re^{i\theta
})\right\Vert _{\mathcal{L}(X,Y)}^{2}dm(\theta ) &\leq &sup_{0\leq r<1}\left[
\dint\limits_{\mathbb{T}}\left[ \dint\limits_{\mathbb{T}}\left\Vert \phi
_{f}(t)\right\Vert _{\mathcal{L}(X,Y)}^{2}P_{re^{i\theta }}(t)dm(t)\right]
dm(\theta )\right]  \\
&=&sup_{0\leq r<1}\left[ \dint\limits_{\mathbb{T}}\left[ \dint\limits_{%
\mathbb{T}}P_{re^{i\theta }}(t)dm(\theta )\right] \left\Vert \phi
_{f}(t)\right\Vert _{\mathcal{L}(X,Y)}^{2}dm(t)\right]  \\
&=&\dint\limits_{\mathbb{T}}\left\Vert \phi _{f}(t)\right\Vert _{\mathcal{L}%
(X,Y)}^{2}dm(t),
\end{eqnarray*}

car $\dint\limits_{\mathbb{T}}P_{re^{i\theta }}(t)dm(\theta )=1.$

Montrons l'in\'{e}galit\'{e} inverse.

Pour presque tout $\theta \in \mathbb{T}$ et tout $x\in B_{X}$ nous avons $%
\phi _{f}(\theta )(x)=lim_{r\rightarrow ^{-}1}f(re^{i\theta })(x).$
Observons maintenant que la fonction $z\in \mathbb{D}\rightarrow \left\Vert
f(z)\right\Vert _{\mathcal{L}(X,Y)}$ est sous-harmonique et que  $%
\sup_{0\leq r<1}\dint\limits_{\mathbb{T}}\left\Vert f(re^{i\theta
})\right\Vert _{\mathcal{L}(X,Y)}^{2}dm(\theta )<+\infty $. Ceci entra\^{\i}%
ne que $lim_{r\rightarrow ^{-}1}\left\Vert f(re^{i\theta })\right\Vert _{%
\mathcal{L}(X,Y)}$ existe pour presque tout $\theta \in \mathbb{T}.$ Donc $%
\left\Vert \phi _{f}(\theta )(x)\right\Vert _{Y}=lim_{r\rightarrow
^{-}1}\left\Vert f(re^{i\theta })(x)\right\Vert _{Y}\leq lim_{r\rightarrow
^{-}1}\left\Vert f(re^{i\theta })\right\Vert _{\mathcal{L}(X,Y)},$ pour tout 
$x\in B_{X}.$

D'apr\`{e}s ce qui est pr\'{e}c\`{e}de, nous avons $\left\Vert \phi
_{f}(\theta )\right\Vert _{\mathcal{L}(X,Y)}=\sup_{x\in B_{X}}\left\Vert
\phi _{f}(\theta )(x)\right\Vert _{Y}\leq lim_{r\rightarrow ^{-}1}\left\Vert
f(re^{i\theta })\right\Vert _{\mathcal{L}(X,Y)}$ pour presque tout $\theta
\in \mathbb{T}.$ Ceci implique d'apr\`{e}s le lemme de Fatou que 
\begin{eqnarray*}
\dint\limits_{\mathbb{T}}\left\Vert \phi _{f}(\theta )\right\Vert _{\mathcal{%
L}(X,Y)}^{2}dm(\theta ) &\leq &\dint\limits_{\mathbb{T}}lim_{r\rightarrow
^{-}1}\left\Vert f(re^{i\theta })\right\Vert _{\mathcal{L}%
(X,Y)}^{2}dm(\theta ) \\
&\leq &\underset{r\rightarrow ^{-}1}{\lim }\dint\limits_{\mathbb{T}%
}\left\Vert f(re^{i\theta })\right\Vert _{\mathcal{L}(X,Y)}^{2}dm(\theta )=
\\
&&sup_{0\leq r<1}\dint\limits_{\mathbb{T}}\left\Vert f(re^{i\theta
})\right\Vert _{\mathcal{L}(X,Y)}^{2}dm(\theta ).
\end{eqnarray*}

Par cons\'{e}quent $\left\Vert f\right\Vert _{\mathbf{h}^{2}(\mathbb{D},%
\mathcal{L}(X,Y))}^{2}=\dint\limits_{\mathbb{T}}\left\Vert \phi _{f}(\theta
)\right\Vert _{\mathcal{L}(X,Y)}^{2}dm(\theta ).\blacksquare $

\emph{Etape 3:}\textbf{\ }Soient $f\in \mathbf{h}^{2}\left[ \mathbb{D},%
\mathcal{L}(X,Y)\right] )$ et $(\rho _{n})_{n\geq 0}$ une suite dans $\left]
0,1\right[ $ convergeant vers 1$.$ Montrons qu'il existe une sous-suite $%
(\rho _{n_{k}})_{k\geq 0}$ telle que pour presque tout $\theta \in \mathbb{T}
$ $Vf(\rho _{n_{k}}e^{i\theta })\underset{k\rightarrow +\infty }{\rightarrow 
}V\phi _{f}(\theta )$ faiblement dans $Z.$

Pour tout $n\in \mathbb{N}$ notons $f_{n}(z)=f(\rho _{n}z),$ $z\in \mathbb{D}%
.$

\bigskip\ Observons que%
\begin{eqnarray*}
\overline{lim}_{n\rightarrow +\infty }\left\Vert \frac{(f_{n}+f)}{2}%
\right\Vert _{\mathbf{h}^{2}(\mathbb{D},\mathcal{L}(X,Y))} &\leq &\frac{1}{2}%
\times \lim_{n\rightarrow +\infty }\left\Vert f_{n}\right\Vert _{\mathbf{h}%
^{2}(\mathbb{D},\mathcal{L}(X,Y))}^{{}}+ \\
\frac{1}{2}\left\Vert f\right\Vert _{\mathbf{h}^{2}(\mathbb{D},\mathcal{L}%
(X,Y))} &=&\left\Vert f\right\Vert _{\mathbf{h}^{2}(\mathbb{D},\mathcal{L}%
(X,Y))}.
\end{eqnarray*}

D'autre part, pour tout $r\in \left] 0,1\right[ $%
\begin{equation*}
\left\Vert \frac{(f_{n}+f)}{2}\right\Vert _{\mathbf{h}^{2}(\mathbb{D},%
\mathcal{L}(X,Y))}^{2}\geq \dint\limits_{\mathbb{T}}\left\Vert \frac{\left[
f(\rho _{n}re^{i\theta })+f(re^{i\theta })\right] }{2}\right\Vert _{\mathcal{%
L}(X,Y)}^{2}dm(\theta ),
\end{equation*}

d'apr\`{e}s le lemme de Fatou, pour tout $r\in \left[ 0,1\right[ $ nous avons

\begin{eqnarray*}
\underline{\lim }_{n\rightarrow +\infty }\left\Vert (\frac{f_{n}+f)}{2}%
\right\Vert _{\mathbf{h}^{2}(\mathbb{D},\mathcal{L}(X,Y))}^{2} &\geq &%
\underline{\lim }_{n\rightarrow +\infty }\dint\limits_{\mathbb{T}}\left\Vert 
\frac{\left[ f(\rho _{n}re^{i\theta })+f(re^{i\theta })\right] }{2}%
\right\Vert _{\mathcal{L}(X,Y)}^{2}dm(\theta ) \\
&\geq &\dint\limits_{\mathbb{T}}\left\Vert f(re^{i\theta })\right\Vert _{%
\mathcal{L}(X,Y)}^{2}dm(\theta ).
\end{eqnarray*}

Il en r\'{e}sulte que

\begin{equation*}
\lim_{n\rightarrow +\infty }\left\Vert \frac{(f_{n}+f)}{2}\right\Vert _{%
\mathbf{h}^{2}(\mathbb{D},\mathcal{L}(X,Y))}^{2}=\left\Vert f\right\Vert _{%
\mathbf{h}^{2}(\mathbb{D},\mathcal{L}(X,Y))}^{2},
\end{equation*}

\bigskip\ nous d\'{e}duisons que

\begin{equation*}
\frac{1}{2}\left\Vert (f_{n}\right\Vert _{\mathbf{h}^{2}(\mathbb{D},\mathcal{%
L}(X,Y))}^{2}+\frac{1}{2}\left\Vert f\right\Vert _{\mathbf{h}^{2}(\mathbb{D},%
\mathcal{L}(X,Y))}^{2}-\left\Vert \frac{(f_{n}+f)}{2}\right\Vert _{\mathbf{h}%
^{2}(\mathbb{D},\mathcal{L}(X,Y))}^{2}\underset{n\rightarrow +\infty }{%
\rightarrow }0,
\end{equation*}

(remarquons que $\left\Vert (f_{n}\right\Vert _{\mathbf{h}^{2}(\mathbb{D},%
\mathcal{L}(X,Y))}^{2}\underset{n\rightarrow \infty }{\rightarrow }%
\left\Vert f\right\Vert _{\mathbf{h}^{2}(\mathbb{D},\mathcal{L}(X,Y))}^{2}).$
D'apr\`{e}s l'\'{e}tape 2, nous obtenons que

\begin{equation*}
\dint\limits_{\mathbb{T}}\left[ \frac{1}{2}\left\Vert f(\rho _{n}e^{i\theta
})\right\Vert _{\mathcal{L}(X,Y)}^{2}+\frac{1}{2}\left\Vert \phi _{f}(\theta
)\right\Vert _{\mathcal{L}(X,Y)}^{2}-\left\Vert \frac{\left[ f(\rho
_{n}e^{i\theta })+\phi _{f}(\theta )\right] }{2}\right\Vert ^{2}\right]
dm(\theta )\underset{n\rightarrow +\infty }{\rightarrow }0.
\end{equation*}

Par cons\'{e}quent il existe une sous-suite $(\rho _{n_{k}})_{k\geq 0}$
telle que%
\begin{equation*}
\frac{1}{2}\left\Vert f(\rho _{n_{k}}e^{i\theta })\right\Vert _{\mathcal{L}%
(X,Y)}^{2}+\frac{1}{2}\left\Vert \phi _{f}(\theta )\right\Vert _{\mathcal{L}%
(X,Y)}^{2}-\left\Vert \frac{\left[ f(\rho _{n_{k}}e^{i\theta })+\phi
_{f}(\theta )\right] }{2}\right\Vert _{\mathcal{L}(X,Y)}^{2}\underset{%
k\rightarrow +\infty }{\rightarrow }0,
\end{equation*}

pour presque tout $\theta \in \mathbb{T}$ (remarquons que $\frac{1}{2}%
\left\Vert f(\rho _{n}e^{i\theta })\right\Vert _{\mathcal{L}(X,Y)}^{2}+\frac{%
1}{2}\left\Vert \phi _{f}(\theta )\right\Vert _{\mathcal{L}%
(X,Y)}^{2}-\left\Vert \frac{\left[ f(\rho _{n}e^{i\theta })+\phi _{f}(\theta
)\right] }{2}\right\Vert _{\mathcal{L}(X,Y)}^{2}\geq 0,$ pour presque tout $%
\theta \in \mathbb{T}$ et tout $n\in \mathbb{N)}.$ Comme $V$ est faiblement $%
LUR$, pour presque tout $\theta \in \mathbb{T}$ nous avons $Vf(\rho
_{n_{k}}e^{i\theta })\underset{k\rightarrow +\infty }{\rightarrow }V\phi
_{f}(\theta )$ faiblement dans $Z$ pour presque tout $\theta \in \mathbb{T}.$

\emph{Etape 4:} Montrons que $Vf$ $(re^{i\theta })\underset{r\rightarrow
^{-}1}{\rightarrow }V\phi _{f}(\theta )$ dans $Z$ pour presque tout $\theta
\in \mathbb{T}.$

D'apr\`{e}s le th\'{e}or\`{e}me de mesurabilit\'{e} de Pettis \cite[th. II 2]%
{Du}, la fonction $V\phi _{f}$ est $m-mesurable$ pour la topologie de la
norme de $Z.$ Donc, pour tout $u^{\ast }\in Z^{\ast }$ et tout $z\in \mathbb{%
D}$ $\dint\limits_{\mathbb{T}}(V\phi _{f}(t),u^{\ast
})P_{z}(t)dm(t)=(Vf(z),u^{\ast })$, car pour presque tout $\theta \in 
\mathbb{T}$ $Vf(\rho _{n_{k}}e^{i\theta })\underset{k\rightarrow +\infty }{%
\rightarrow }V\phi _{f}(\theta )$ faiblement dans $Z,$ donc $\dint\limits_{%
\mathbb{T}}V\phi _{f}(t)P_{z}(t)dm(t)=Vf(z),$ pour tout $z\in \mathbb{D}.$
D'apr\`{e}s \cite[th.7.6]{Ru} $Vf$ $(re^{i\theta })\underset{r\rightarrow
^{-}1}{\rightarrow }V\phi _{f}(\theta )$ pour presque tout $\theta \in 
\mathbb{T}.\blacksquare $

\begin{corollary}
\label{nb}Supposons que $X$ soit s\'{e}parable, $Y$ a la propri\'{e}t\'{e}
de Radon-Nikodym et $\mathcal{L}(X,Y)$ est faiblement $LUR.$ Alors $\mathcal{%
L}(X,Y)$ a la propri\'{e}t\'{e} de Radon-Nikodym.
\end{corollary}

\begin{theorem}
\label{cs}Soient $Z$ un espace de Banach et $V:X\rightarrow Y$ un op\'{e}%
rateur born\'{e}. On d\'{e}finit l'op\'{e}rateur $\widetilde{V}:\mathcal{L}%
(Y,Z)\rightarrow \mathcal{L}(X,Z)$ par $\widetilde{V}T=T\circ V,$ $T\in 
\mathcal{L}(Y,Z).$ Supposons que $\widetilde{V}$ soit \`{a} valeurs dans un
sous-espace $WCG$ $E$ de $\mathcal{L}(X,Z),$ $Z$ a la propri\'{e}t\'{e} de
de Radon-Nikodym et que $Y$ soit s\'{e}parable. Alors $\widetilde{V}$ est un
op\'{e}rateur de Radon-Nikodym.
\end{theorem}

D\'{e}montration.

Soit $f\in \mathbf{h}^{2}\left[ \mathbb{D},\mathcal{L}(Y,Z)\right] .$ Par un
argument analogue \`{a} celui du th\'{e}or\`{e}me \ref{oq}, on montre qu'il
existe une fonction $\phi _{f}$ \`{a} valeurs presque partout dans $\mathcal{%
L}(Y,Z)$ telle que pour presque tout $\theta \in \mathbb{T}$ et tout $y\in Y$
$f(re^{i\theta })y\underset{r\rightarrow ^{-}1}{\rightarrow }\phi
_{f}(\theta )y$ dans $Z.$ Ceci entra\^{\i}ne que pour presque tout $\theta
\in \mathbb{T}$ et tout $x\in X$ $f(re^{i\theta })Vx\underset{r\rightarrow
^{-}1}{\rightarrow }\phi _{f}(\theta )Vx$ dans $Z.$

\emph{Etape 1:}\textbf{\ }Montrons que $\phi _{f}(.)\circ V$ est faiblement
mesurable dans $E$.

Consid\'{e}rons%
\begin{equation*}
L=\left\{ u^{\ast }\in E^{\ast };\text{ l'application }(\phi _{f}(.)\circ
V,u^{\ast })\text{ est }m-mesurable\right\} .
\end{equation*}

Soit $U\in E.$ Supposons que pour tout $x\otimes v^{\ast }\in X\otimes
Z^{\ast }$ $(U,x\otimes z^{\ast })=(U(x),v^{\ast })=0.$ Ceci impique que $%
U=0.$ Donc $X\otimes Z^{\ast }$ est pr\'{e}faiblement dense dans $E^{\ast }.$

On remarque d'autre part, que $X\otimes Z^{\ast }\subset L,$ il suffit donc
montrer que $L$ est pr\'{e}faiblement ferm\'{e} dans $E^{\ast }$,.D'apr\`{e}%
s le th\'{e}or\`{e}me de Banach-Dieudonn\'{e}, il suffit de montrer que la
boule unit\'{e} de $L$ est pr\'{e}faiblement ferm\'{e}e dans $E^{\ast }.$

\bigskip Comme $E$ est un espace $WCG,$ il existe un espace de Banach $F$ r%
\'{e}flexif et $S:F\rightarrow E$ un op\'{e}rateur born\'{e} d'image dense.
Soit maintenant ($u_{\alpha }^{\ast })_{\alpha \in I}$ une suite g\'{e}n\'{e}%
ralis\'{e}e dans la boule unit\'{e} de $L$ telle que $u_{\alpha }^{\ast
}\rightarrow u^{\ast }\in E^{\ast }$ pr\'{e}faiblement dans $E^{\ast }.$
Cela implique que $S^{\ast }u_{\alpha }^{\ast }\rightarrow S^{\ast }u^{\ast }
$ faiblement dans $F^{\ast }.$ Il existe alors une suite $(t_{n}^{\ast
})_{n\geq 0}$ dans l'enveloppe convexe de $(u_{\alpha }^{\ast })_{a\in I}$
telle que $S^{\ast }t_{n}^{\ast }\underset{n\rightarrow +\infty }{%
\rightarrow }S^{\ast }u^{\ast }$ fortement dans $F^{\ast },$ ceci entr\^{\i}%
ne que $t_{n}^{\ast }\underset{n\rightarrow +\infty }{\rightarrow }u^{\ast }$
$\sigma (E^{\ast },S(F)).$ Donc pr\'{e}faiblement dans $E^{\ast }.$
Observons que $L$ est pr\'{e}faiblement s\'{e}quentiellement ferm\'{e} dans $%
E^{\ast }$, il en r\'{e}sulte que $u^{\ast }\in L.$

Comme $E$ est un espace $WCG$, d'apr\`{e}s le r\'{e}sultat de D.R Lewis \cite%
[p.88]{Du} il existe une fonction $\psi $ qui est $m-mesurable$ \`{a}
valeurs dans $E$, telle que pour tout $u^{\ast }\in E^{\ast },$ en
particulier tout $x\otimes v^{\ast }\in X\otimes Z^{\ast }$ $(\phi
_{f}(.)V(x),v^{\ast })=(\psi (.)(x),v^{\ast })$ p.s.

\emph{Etape 2}\textbf{: }Montrons que $\dint\limits_{\mathbb{T}}\psi
(t)P_{z}(t)dm(t)=f(z)\circ V,$ $z\in \mathbb{D}.$

Pour tout $x\otimes v^{\ast }\in X\otimes Z^{\ast }$

\begin{eqnarray*}
(\dint\limits_{\mathbb{T}}\psi (t)P_{z}(t)dm(t),x\otimes v^{\ast })
&=&\dint\limits_{\mathbb{T}}(\psi (t)x,v^{\ast })P_{z}(t)dm(t)) \\
&=&\dint\limits_{\mathbb{T}}(\phi _{f}(t)\circ V(x),v^{\ast
})P_{z}(t)dm(t))=(f(z)\circ V(x),v^{\ast }),\text{ }z\in \mathbb{D}.
\end{eqnarray*}

Par cons\'{e}quent $\dint\limits_{\mathbb{T}}\psi (t)P_{z}(t)dm(t)=f(z)\circ
V$, \ $z\in \mathbb{D}.\blacksquare $

En choisissant $X=Y$ et $T$ l'identit\'{e} de $X,$ on tire le corollaire
suivant:

\begin{corollary}
\label{ghm}Soit $Z$ un espace de Banach. Supposons que $X$ soit un espace de
Banach s\'{e}parable, $Z$ a la propri\'{e}t\'{e} de Radon-Nikodym et que $%
\mathcal{L}(X,Z)$ soit un espace $WCG.$ Alors $\mathcal{L}(X,Z)$ a la propri%
\'{e}t\'{e} de Radon-Nikodym.
\end{corollary}

\bigskip

\end{document}